\documentclass[11pt,twoside]{article}
\usepackage{latexsym}
\usepackage{amssymb,amsbsy,amsmath,amsfonts,amssymb,amscd}
\usepackage{epsfig, graphicx}
\usepackage{color}
\newcommand  \ind[1]  {   {1\hspace{-1.2mm}{\rm I}}_{\{#1\} }    }
\newcommand{\fer}[1]{(\ref{#1})}

\newcommand{\commentout}[1]{}
\newcommand {\al} {\alpha}

\newcommand {\sg} {\sigma}

\newcommand {\Chi} {{\bf \raise 2pt \hbox{$\chi$}} }
\newcommand {\f}   {\frac}
\newcommand {\p}   {\partial}
\newcommand{\dis}{\displaystyle}

\newcommand{\beq}{\begin{equation}}
\newcommand{\beqa}{\begin{array}{l}}
\newcommand{\bea} {\begin{array}{rl}}
\newcommand{\eeq}{\end{equation}}
\newcommand{\eeqa}{\end{array}}
\newcommand{\eea} {\end{array}}
\newtheorem{theorem}{Theorem}[section]

\newcommand{\qed}{{ \hfill
                       {\unskip\kern 6pt\penalty 500
                       \raise -2pt\hbox{\vrule\vbox to 6pt{\hrule width 6pt
                       \vfill\hrule}\vrule} \par}   }}
\title{\Large \bf Exponential decay for the growth-fragmentation/cell-division equation}

\author{Philippe Lauren\c cot  \thanks{Institut de
Math\'ematiques de Toulouse, CNRS UMR~5219, Universit\'e de Toulouse, F--31062 Toulouse cedex 9, France. 
E-mail: laurenco@math.univ-toulouse.fr}
\and Beno\^ \i t Perthame\thanks{
Universit\'e Pierre et Marie Curie-Paris 6, UMR 7598 LJLL, BC187, 4, place Jussieu,  F-75252 Paris cedex 5, and Institut Universitaire de France. E-mail: benoit.perthame@upmc.fr}
}

\date{\today}

\begin{document}
\maketitle
\pagestyle{plain}
\pagenumbering{arabic}

\begin{abstract}
We consider the linear growth-fragmentation equation arising in the modelling of cell division or polymerisation processes. For constant coefficients, we prove that the dynamics converges to the steady state with an exponential rate. The control on the initial data uses an elaborate $L^1$-norm that seems to be necessary. It also reflects the main idea of the proof which is to use an anti-derivative of the solution. The main technical difficulty is related to the  entropy dissipation rate which is too weak to produce a Poincar\'e inequality.  
\end{abstract}

\medskip

\noindent {\bf Keywords} Growth-fragmentation equations, cell division, polymerisation process, size repartition, convergence to equilibrium, temporal rate of convergence.
\\[5mm]
{\bf AMS Class. No.} 35B40, 45K05, 82D60, 92D25

\section{Introduction}
\label{sec:intro}

We consider the growth-fragmentation (or cell division) equation

\beq \label{eq:gf}
\left\{ \beqa
\displaystyle{\f \p {\p t} n(x,t) + \f \p {\p x} n(x,t) + kB n(x,t) = B\int_0^\infty \kappa(x,y) n(y,t) dy },
\\[4mm]
n(x=0,t) = 0,
\eeqa\right.
\eeq
together with an initial condition $n(x, t=0)=n^0(x)\in L^1(0,\infty)$. Here, the number $k>1$ represents the average number of pieces after fragmentation, while $B>0$ measures the relative intensity of the fragmentation process versus the growth process and $x \mapsto \kappa(x,y)$ gives the size repartition after fragmentation of a polymer (or division of a cell) of size $y$. Throughout this paper, the fragment distribution function $\kappa$ is assumed to enjoy the following properties:
\beq \label{as1}
 \kappa(x,y)= 0 \; \text{ for } 0<y < x,
\eeq
\beq \label{as2}
 \int_{0}^{y} \kappa(x,y)dx=k , \qquad   \int_{0}^{y} x \kappa(x,y)dx=  y .
\eeq
These assumptions are usual in modelling fragmentation phenomena \cite{BP} and lead to two balance laws for the number of fragments $\int_0^\infty n(x,t)dx$ and the total mass $\int_0^\infty x \; n(x,t)dx$, namely 
\beq \label{eq:number}
\int_0^\infty n(x,t) dx=\varrho := \int_0^\infty n^0(x)dx \qquad \mbox{ for all }\;\; t \geq 0, 
\eeq
\beq \label{eq:mass}
\f{d}{dt} \int_0^\infty x \; n(x,t)dx+(k-1) B \int_0^\infty x \; n(x,t)dx = \int_0^\infty n(x,t)dx = \varrho \qquad \mbox{ for all }\;\; t \geq 0.
\eeq
These are instrumental for the existence theory of weak solutions.

The evolution equation \fer{eq:gf} is already renormalised so that zero is the first eigenvalue of the underlying linear integro-differential operator. Indeed, recall that, under fairly general assumptions on $\kappa$ that were stated for instance in \cite{M,BP}, there exists a unique steady state $N\in L^1(0,\infty)$ satisfying
\beq \label{eq:gfe}
\left\{ \beqa
\displaystyle{\f \p {\p x} N(x) + k B N(x) =B \int_0^\infty \kappa(x,y) N(y) dy ,}
\\[4mm]
\displaystyle{N(x=0) = 0, \qquad N>0 , \quad \int_0^\infty N(x) dx=1.}
\eeqa\right.
\eeq
In addition, convergence of solutions $n$ to \fer{eq:gf} with initial number of fragments $\varrho$ towards $\varrho N$ is established in \cite[Section~4]{MMP} under suitable assumptions on $\kappa$ but without providing a rate of convergence (see also \cite[Chapter~4]{BP}). These convergence results are based on the analysis of the entropy dissipation  that we recall in Section~\ref{sec:examples}. However, an estimate on the time needed to reach the steady state is of importance: in particular, an exponential rate of convergence guarantees that the steady state is reached fast enough to be observed in practice. Another motivation to investigate the rate of convergence is the study of the stability of steady states in nonlinear cases, at least in weakly nonlinear regimes. Such nonlinear models are met for example in prion dynamics \cite{CLODMLP,GPW2006,LW}, cell division with finite resources \cite{MD,BP}, network communications \cite{baccelli}, and many other fields.

There are two particular  cases where an exponential decay is proved. The first is equal mitosis equations (each cell divides in two cells of the same size) which corresponds to the choice $\kappa(x,y)= 4 \delta(y-2x)$ (see Section~\ref{sec:examples}), and the method of proof in \cite{BP,PR} is extended here. The second is the case of age structured equations which is very particular and more classical because it corresponds to usual modelling; it corresponds to the formal limit $\sg \to 0$ in the case of general mitosis below in Section~\ref{sec:examples} and the exponential decay rate is known even the case of non-constant coefficients, see \cite{feller,GwP,MD}.

The aim of this note is to show that such an exponential decay is actually valid for a broader class of fragment distribution functions $\kappa$, namely those satisfying the positivity condition
\beq \label{as3}
\beta(x,y):= - \f \p {\p y}  \int_0^x \kappa(z,y)dz \geq 0,
\eeq
in addition to \fer{as1} and \fer{as2}. Below we give several examples for which these properties are met. The interest of our result is not only to extend the exponential trend to equilibrium to a wider class of data but also to explain the general structure by the positivity property in \fer{as3}. The semi-norm on the initial data that comes out from our analysis is defined by
\begin{equation}
||| n|||  :=  \int_0^\infty \left|  \frac{\p M}{\p x}(x) + kB M(x) -B \int_0^\infty \beta(x,y) M(y)dy \right| dx+ (k+1)B \int_0^\infty |M(x) | dx  
\label{eq:norm}
\end{equation}
with $M(x)=\int_0^x (n(z)-\varrho \; N(z)) dz$ and $\varrho= \int_0^\infty n(y) dy$ (following the notation introduced in  \fer{eq:number}).

\begin{theorem} With assumptions \fer{as1}, \fer{as2}, \fer{as3} and an initial condition satisfying $||| n^0||| < \infty$, the corresponding solution $n$ to \fer{eq:gf} satisfies\beq \label{eq:decayn}
\int_0^\infty |n(x,t) -\varrho N(x) |dx \leq ||| n^0||| \; e^{-(k-1)Bt} \qquad \mbox{ for all }\;\; t \geq 0.
\eeq
\label{th:decay}
\end{theorem}

Up to our knowledge it is not possible to obtain exponential decay with only control of $L^p$-norms on the initial data (even though we are not aware of an explicit counterexample). This is related to the main difficulty in analysing the decay rate; the entropy dissipation is too weak to produce a Poincar\'e inequality. We explain this in more details in Section \ref{sec:examples}. This is usual and the case of Fokker-Planck equations and ``hypocoercive'' equations are other examples \cite{DV,Villani}.

The outline of the paper is as follows. In the next section, we give some examples where the assumptions \fer{as1}, \fer{as2}, \fer{as3} are satisfied. Then we prove Theorem~\ref{th:decay}. The proof uses an integral change of unknown which has already been used in \cite{L05,BP,PR} for similar purposes (see also \cite{FL06,LW,NP05} where it is applied to the uniqueness issue for related models). More precisely, we study the equation solved by an anti-derivative of the solution to \fer{eq:gf}. A crucial point here is that it satisfies a closed equation, a property due to the constant fragmentation coefficient $kB$. We gather several remarks in the last section.

\section{Examples and entropy dissipation}
\label{sec:examples}

We present here several standard examples of coefficients $\kappa(x,y)$ that our method is able to handle. We also recall the relative entropy associated with the equations at hand and indicate why a Poincar\'e inequality cannot hold true.	
\\
\\
{\bf Equal mitosis.} The most general model describing cell division relies on the choice  
\beq \label{eq:ex1}
  \kappa(x,y)= 4 \delta(y-2x)= 2  \delta(x-y/2), \qquad k=2,
\eeq
that is, each cell divides in two cells of equal sizes. Assumptions \fer{as1} and \fer{as2} are obviously satisfied. Furthermore, we can compute
$$
 \int_0^x \kappa(z,y)dz = 2 \ind{x\geq y/2}= 2 \ind{y \leq 2x}, \qquad 
\beta(x,y) = 2  \delta(y-2x)=   \delta(x- y/2),
$$ 
and \fer{as3} is satisfied. 
\\[5mm]
{\bf Uniform fragmentation.} This is the simplest and most standard model for (physical or biological) polymers. It assumes uniform distribution of the fragments after each fragmentation event, i.e.,   
\beq \label{eq:ex2}
 \kappa(x,y)= \f 2 y \ind{x\leq y}, \qquad k=2.
\eeq
Assumptions \fer{as1} and \fer{as2} are again obviously satisfied. In addition,
$$
 \int_0^x \kappa(z,y)dz = 2\  \f {\min\{x,y\}}{y} , \qquad
\beta(x,y) =  2\ \f x {y^2} \ind{y>x} ,
$$
and \fer{as3} is also satisfied. 
\\[5mm]
{\bf General Mitosis.}  An extension of equal mitosis is to take, for some parameter $\sigma$ with $0 < \sigma \leq 1/2$,  
\beq \label{eq:ex3}
\kappa(x,y)=  \f 1 {\sigma} \delta\left( y- \f x \sigma \right)+ \f 1 {1-\sigma} \delta\left( y-\f x {1-\sigma} \right)=   \delta(x-\sigma y )+ \delta(x- (1-\sigma)y).
\eeq
Assumptions \fer{as1} and \fer{as2} are clearly fulfilled, still with $k=2$, and
$$
\int_0^x \kappa(z,y)dz =  \ind{x\geq \sigma y} +  \ind{x\geq (1-\sigma) y} = \ind{y \leq x/\sigma } +  \ind{y \leq x/(1-\sigma)}, 
$$
$$
\beta(x,y) =   \delta\left( y- \f{x}{\sigma} \right) +  \delta\left( y- \f{x}{1-\sigma} \right)=  \sigma \delta(x- \sigma y)+ (1-\sigma)\delta(x- (1- \sigma) y).
$$
Therefore \fer{as3} is again satisfied. 
\\[5mm]
{\bf Homogeneous fragmentation.} This generalisation of uniform fragmentation consists in using 
\beq \label{eq:ex4}
\kappa(x,y)= (2+\al)\;  \f {x^\al}{ y^{1+\al}} \ind{x\leq y}, \qquad k=\f{2+\al}{1+\al}.
\eeq
for some $\al>-1$ and satisfies \fer{as1} and \fer{as2}. Also
$$
\int_0^x \kappa(z,y)dz = k  \f {\min\{x,y\}^{1+\al}} {y^{1+\al}}  , \qquad \beta(x,y) =(2+\al) \f {x^{1+\al}} {y^{2+\al}}  {\ind{y>x}}  ,
$$
and \fer{as3} is again fulfilled. 
\\

All these examples illustrate the main difficulty in order to prove the exponential decay of solutions and why entropy methods are not enough. We recall from \cite{MMP,BP, Villani} and the references therein that the relative entropy inequality for equation \fer{eq:gf} reads
$$
\f{d}{dt}\int_0^\infty N(x) H\left(\f{n(x,t)}{N(x)}\right) dx = - D_H[n(t)]
$$
$$
D_H[n]:= B \int_0^\infty \int_0^\infty \kappa(x,y) N(y) \left[
H\left(\f{n(y)}{N(y)}\right) -H\left(\f{n(x)}{N(x)}\right)-H'\left(\f{n(x)}{N(x)}\right) \left(\f{n(y)}{N(y)}-\f{n(x)}{N(x)}\right)
\right] dx \; dy.
$$

For a convex function $H$ this is a non-negative quantity and this entropy inequality is usually enough to prove the strong convergence of $n(t)$ to $\varrho N$ in a norm that depends on the choice of $H$. However the simplest method to prove exponential decay is through a Poincar\'e inequality (whenever available) which asserts the existence of $\nu>0$ such that
$$
\nu \int_0^\infty N(x) H\left(\f{n(x)}{N(x)}\right) dx \leq D_H[n],
$$
for all functions $n$ satisfying $\int_0^\infty n(x) dx=0$.

The usual setting for the quadratic entropy $H(r)=r^2$ is to say that $u = n/N$ satisfies
$$
\nu \int_0^\infty N(x) u(x)^2 dx \leq  D_2[u]:= B \int_0^\infty \int_0^\infty \kappa(x,y) N(y)  \big( u(x)-u(y)\big)^2 dx \; dy, 
$$
whenever $\int_0^\infty u(x) N(x) dx=0$. As far as we know, the validity of such an inequality is unclear and it is actually false in the case of equal mitosis. Indeed, in that case, $\kappa$ is given by \eqref{eq:ex1} and the right-hand side $D_2[u]$ is
$$
D_2[u]= 4B \int_0^\infty N(2x) \big( u(x)-u(2x)\big)^2 dx,
$$
which clearly vanishes for the class of functions that satisfy $u(x)=u(2 x)$ for all $x\in (0,\infty)$: such functions are obtained from a given bounded function $\phi$ defined on $[1,2]$ by setting $u(x)=\phi(x/2^j)$ for $x\in [2^j,2^{j+1})$ and $j\in \mathbb{Z}$.

\section{Proof of Theorem~\ref{th:decay}}
\label{sec:proof}

As already mentioned, the proof uses an integral change of unknown which is our starting point.

\noindent {\em Step 1. (change of unknown function)}
Let $n$ be a solution to \fer{eq:gf} and set
$$
M(x,t) = \int_0^x [n(z,t) -\varrho N(z)] dz \;\;\mbox{ for }\;\; (x,t)\in [0,\infty)^{2}.
$$
It satisfies $M(0,t)=M(\infty, t)=0$ thanks to the conservation of the number of fragments \fer{eq:number} which implies  
$$
\int_0^\infty  [n(z,t) -\varrho N(z)] dz= \int_0^\infty  n^0(z)dz -\rho \int_0^\infty N(z) dz=0.
$$

The motivation for using $M$ is that it satisfies a closed equation of the same type, namely,
\beq \label{eq:m}
\left\{ \beqa
\displaystyle{\f \p {\p t} M(x,t) + \f \p {\p x} M(x,t) + kB M(x,t) = B  \dis \int_0^\infty \beta(x,y)  M(y,t) dy,}
\\[4mm]
M(x=0, t) = 0.
\eeqa\right.
\eeq
Indeed, integrating the first equation in \fer{eq:gf} over $(0,x)$, $x\in (0,\infty)$, we obtain 
\begin{eqnarray*}
\f \p {\p t} M(x,t) + \f \p {\p x} M(x,t) + kB M(x,t) &=& B\int_{z=0}^x \int_{y=0}^\infty \kappa(z,y) n(y,t) dy dz \\
&=&B \dis   \int_{y=0}^\infty \f \p {\p y} M(y,t) \int_{z=0}^{x} \kappa(z,y) dz dy \\
&=& B  \dis \int_{y=0}^\infty \left[- \f \p {\p y} \int_{z=0}^{x}  \kappa(z,y) dz \right] M(y,t) dy  \\ 
&=& B  \dis \int_0^\infty \beta(x,y)  M(y,t) dy .
\end{eqnarray*}

\medskip

\noindent {\em Step 2. (Property of $\beta$)}
To proceed further, we need a specific property of $\beta$, namely
\begin{eqnarray}
\int_0^\infty \beta(x,y) dx&=& - \int_{x=0}^{y} \f \p {\p y} \left( \int_{z=0}^{x}   \kappa(z,y) dz \right) dx \nonumber\\
&=& - \f \p {\p y} \left( \int_{x=0}^y  \int_{z=0}^{x}  \kappa(z,y) dzdx \right) + \int_{z=0}^{y} \kappa(z,y) dz \nonumber\\
&=& - \f \p {\p y} \left( \int_{z=0}^y \kappa(z,y)  \int_{x=z}^{y}  dx dz \right) + k\nonumber\\
&=& - \f \p {\p y} \left( \int_{z=0}^y  (y-z)   \kappa(z,y) dz \right) + k \nonumber\\
&=& - \f \p {\p y} [(k-1)y ] +k \nonumber
\end{eqnarray}
and thus
\begin{equation}
\int_0^\infty \beta(x,y) dx= 1 . 
\label{pr:beta}
\end{equation}

\medskip

\noindent {\em Step 3. (Exponential decay of $M$)}  
We now argue on the function $M$. Owing to the positivity \fer{as3} of $\beta$, we infer from \fer{eq:m} that 
$$
\f \p {\p t} |M(x,t) |+ \f \p {\p x} |M(x,t) |+ kB |M(x,t)| \leq   \dis \int_{0}^\infty \beta(x,y) | M(y,t) | dy .
$$
After integration of the above inequality with respect to $x$ over $(0,\infty)$, we deduce from \fer{pr:beta} that
\begin{eqnarray*}
\f d{dt} \int_0^\infty  |M(x,t) | dx + kB \int_0^\infty  |M(x,t) | &\leq&   \dis \int_0^\infty \int_{0}^\infty \beta(x,y) | M(y,t) | dydx \\
&\leq& B \int_0^\infty  |M(x,t) | dx ,
\end{eqnarray*}
and thus 
$$
\int_0^\infty  |M(x,t) | dx \leq e^{-(k-1)Bt}\; \int_0^\infty  |M(x,0) | dx \;\;\mbox{ for }\;\; t\ge 0.
$$

\noindent {\em Step 4. (Exponential decay of $n$)}
To transfer the exponential decay to the solution $n$ itself is not always possible, see for instance the coagulation case \cite{L05}. Here, we can follow \cite{PR} and argue as follows. We first notice that ${\p M}/{\p t} $ satisfies the same equation as $M$. Therefore Step~3 and \fer{eq:m} also give
\begin{eqnarray*}
\int_0^\infty  \left| \f{\p M}{\p t}(x,t) \right| dx & \leq & e^{-(k-1)Bt} \; \int_0^\infty  \left| \f{\p M}{\p t}(x,0) \right| dx \\
& \leq & e^{-(k-1)Bt} \; \int_0^\infty\left| \f{\p M}{\p x}(x,0)+ kB M(x,0)-B \int_0^\infty \beta(x,y) M(y,0)dy \right|  dx .
\end{eqnarray*}
Next, owing to \fer{as3} and \fer{pr:beta}, we have
\beq
\label{pr:b2eta}
\int_0^\infty \left| \int_0^\infty \beta(x,y) M(y,t)dy \right| dx \leq \int_0^\infty |M(y,t)| \int_0^\infty \beta(x,y) dx dy = \int_0^\infty |M(y,t)| dy.
\eeq
As a consequence, \fer{eq:m} also reads 
$$
n(x,t)-\varrho N(x)=- \f{\p M}{\p t}(x,t)-kBM(x,t) + \int_0^\infty \beta(x,y) M(y,t)dy,
$$
and we conclude with the help of \fer{pr:b2eta} that 
\begin{eqnarray*}
\int_0^\infty |n(x,t)-\varrho N(x)| dx & \leq & \int_0^\infty  \left| \f{\p M}{\p t}(x,t) \right| dx + (k+1)B  \int_0^\infty  |M(x,t)| dx \\
& \leq & e^{-(k-1)Bt} \; \int_0^\infty\left| \f{\p M}{\p x}(x,0)+ kB M(x,0)-B \int_0^\infty \beta(x,y) M(y,0)dy \right| dx \\
& & + (k+1)B e^{-(k-1)Bt} \;  \int_0^\infty |M(x,0) | dx 
\end{eqnarray*}
as claimed. \qed

\section{Concluding remarks}
\label{sec:cr}

\begin{itemize}
\item[(a)] Owing to \fer{pr:b2eta} and its definition in \fer{eq:norm}, the semi-norm actually satisfies
$$
||| n^{0}||| \le \int_0^\infty\left| \f{\p M}{\p x}(x,0) \right| dx + 2(k+1)B \int_0^\infty | M(x,0)| dx ,
$$
and $|||n^0|||$ is thus dominated by the $W^{1,1}$-norm of $M(.,0)$.
\item[(b)] Theorem~\ref{th:decay} can be extended to growth-fragmentation equations with non-constant velocity
$$
\left\{ \beqa
\displaystyle{\f \p {\p t} n(x,t) + \f \p {\p x} \left( \tau(x) n(x,t) \right) + kB n(x,t) = B\int_0^\infty \kappa(x,y) n(y,t) dy },
\\[4mm]
n(x=0,t) = 0.
\eeqa\right.
$$
The fragment distribution function $\kappa$ is still assumed to fulfil the properties \fer{as1}, \fer{as2}, \fer{as3}, and the velocity $\tau$ is required to satisfy 
$$
\tau(x)\ge \tau_{m}>0 \qquad \vartheta:= (k-1)B - \sup_{x\in (0,\infty)}  \f{\p \tau}{\p x}(x) >0.
$$ 
Under these assumptions, the exponential decay rate is $e^{- \vartheta t}$. Indeed, the equation on $M$ is 
$$
\left\{ \beqa
\displaystyle{\f \p {\p t} M(x,t) +\tau(x) \f \p {\p x} M(x,t) + kB M(x,t) = B  \dis \int_0^\infty \beta(x,y)  M(y,t) dy,}
\\[4mm]
M(x=0, t) = 0.
\eeqa\right.
$$
Therefore, proceeding as in the third step of the proof of Theorem~\ref{th:decay}, we have
$$
\f{d}{dt}\int_0^\infty |M(x,t)|dx +(k-1) B \int_0^\infty |M(x,t)|dx = \int_0^\infty |M(x,t)|  \f{\p \tau}{\p x}(x)dx.
$$
The decay rate for $M$ follows directly from this equality. The positive bound from below required on $\tau$ is next used to transfer this decay rate on $n$, arguing as in the fourth step of the proof of Theorem~\ref{th:decay}. 
\item[(c)] The exponential decay rate in Theorem~\ref{th:decay} seems to be optimal: indeed, consider the case where $\kappa$ is given by \fer{eq:ex1} (equal mitosis), that is, $\kappa(x,y)= 4 \delta(y-2x)= 2  \delta(x-y/2)$, and introduce the function $\xi$ defined by 
$$
\xi(x):= \f{1}{2} \f {\p N}{\p x}\left( \f{x}{2} \right) , \qquad x\in (0,\infty) .
$$
Owing to the properties of $N$ given in \fer{eq:gfe} and the specific choice of $\kappa$, it is straightforward to check that $(x,t)\longmapsto e^{-Bt} \xi(x)$ is the solution to \fer{eq:gf} with initial condition $\xi$ and the integral of $\xi$ vanishes. Consequently, given $\varrho>0$, $n(x,t)=\varrho N(x) + e^{-Bt} \xi(x)$ is the solution to \fer{eq:gf} with initial condition $n^{0}=\varrho N+\xi$ and 
$$
\int_{0}^{\infty} \left| n(x,t) - \varrho N(x) \right| dx = e^{-Bt} \int_{0}^{\infty} |\xi(x)| dx
$$
(recall that $k=2$ in that case).
\item[(d)] Similarly, when $\kappa(x,y)=2 \ind{x\leq y}/y$ (which corresponds to uniform fragmentation, see \eqref{eq:ex2}), the equilibrium $N$ is known explicitly and given by $N(x)=4B^2 x e^{-2Bx}$ for $x\ge 0$. Introducing $\eta(x):=x(Bx-2) e^{-Bx}$ for $x\ge 0$, one readily checks that $(x,t)\longmapsto e^{-Bt} \eta(x)$ is the solution to \fer{eq:gf} with initial condition $\eta$. Observe that the integral of $\eta$ vanishes. As in the previous case, given $\varrho>0$, $n(x,t)=\varrho N(x) + e^{-Bt} \eta(x)$ is the solution to \fer{eq:gf} with initial condition $n^{0}=\varrho N+\eta$. Since 
$$
\int_{0}^{\infty} \left| n(x,t) - \varrho N(x) \right| dx = e^{-Bt} \int_{0}^{\infty} |\eta(x)| dx,
$$
and $k=2$, this example also yields the optimality of the decay rate obtained in Theorem~\ref{th:decay} for that particular case.
\end{itemize}

%
%
%


\end{document}